\documentclass[a4paper,11pt]{article}
\usepackage{amsmath,amssymb,amsfonts}

\newtheorem{theorem}{Theorem}

\newcommand{\N}{\mathbb{N}}
\newcommand{\biz}{\textbf{Proof.\ }}
\newcommand{\kesz}{\hfill$\square$\vspace{0,5cm}}
\renewcommand{\r}{\right)}
\renewcommand{\l}{\left(}
\newcommand{\E}{\mathbb{E}}
\newcommand{\lk}{\left\lbrace}
\newcommand{\rk}{\right\rbrace}
\newcommand{\Z}{\mathbb{Z}}
\renewcommand{\P}{\mathbb{P}}
\newcommand{\tv}{\rightarrow \infty}

\begin{document} 
\thispagestyle{empty}
\begin{center}
\textbf{\begin{Large}Limit distribution of degrees in random family trees                                                                              \end{Large}}\\\vspace{0,3cm}
\textsc{\'Agnes Backhausz}\\
Department of Probability Theory and Statistics\\
E\"{o}tv\" os Lor\'and University\\
P\'azm\'any P\'eter s\'et\'any 1/c, H-1117 Budapest, Hungary \\
Email: \texttt{agnes@cs.elte.hu}
                                                         \end{center}

\begin{abstract}
In a one-parameter model for evolution of random trees, which also includes the Barab\'asi--Albert random tree \cite{BA}, almost sure behavior and the limiting distribution of the degree of a vertex in a fixed position are examined. Results about P\'olya urn models are applied in the proofs.  
\end{abstract}

\noindent Keywords: preferential attachment, random trees, urn models.\\
AMS: 05C80, 60C05, 60F15.

\section{Introduction}

Evolving random graphs and random trees have been widely examined recently, see e.g. \cite{BA, Drmota, Anna}. 
One of the simplest dynamics is the following. At each step one new vertex is born, and it attaches with one 
edge to one of the old vertices. The probability that a given old vertex is chosen is proportional to a fixed 
linear function of its actual degree. The asymptotic degree distribution is well-known. These trees have the 
so-called scale free property: the proportion of vertices of degree $d$ converges to $c_d$ almost surely as the 
number of vertices goes to infinity, and $c_d\sim c\cdot d^{-\gamma} \l d\tv\r$  with some positive constants $c$ and $\gamma$. 

Instead of the degree distribution, we focus on the degree of a vertex in a given position in the tree, as the 
number of vertices goes to infinity. Fix a vertex, e. g. the root of the tree, or the $j$th child of the 
root, or the $k$th child of the $j$th child of the root, etc. $X_n$ denotes the degree of this vertex after $n$ steps. 
We will see that $n^{-\delta}X_n$ converges to a positive random variable almost surely, with some $\delta>0$, 
and we will have some information on the moments and the structure of this random variable. 
We will describe the distribution of these random variables  for the Albert--Barab\'asi tree, where the probability 
that a given old vertex is chosen is proportional to its actual degree. We will also examine a variant, 
the generalized PORT model, where the number of children is relevant instead of the degree.

\section{Random trees}

\subsection{Notations}
Trees are connected graphs without cycles. We assign a vertex, this is the root of the tree. We consider trees growing at discrete time steps. More precisely, we start from the root, and add one new vertex with one edge at each step. When examining the neighbors of a certain vertex, we keep count of the order they were born. Thus our graph is a \textit{rooted ordered tree} (also known as rooted planar tree or family tree). See for example \cite{Drmota}.  
  
We will use the following commonly known terminology and notation for rooted ordered trees \cite{Drmota, Anna}. The vertices are individuals, and the edges of the tree represent the parent-child relations. When a new vertex with one edge is added to the graph, we say that it becomes the child of its only neighbor, its parent. 

We will label the vertices with sequences of positive integers, based on the parent-child relations. We set
\[\Z_+=\lk 1, 2, \ldots\rk, \ \ \ \Z_+^0=\lk \emptyset\rk, \ \ \ \mathcal N=\bigcup_{n=0}^{\infty} \Z_+^n.\]
The label of the root is $\emptyset$. The $j$th children of the root is labelled with $j$. Similarly, the $j$th children of the vertex labelled with $x=\l x_1, \ldots, x_k\r\in \mathcal N$ is labelled with $\l x_1, \ldots, x_k, j\r$. 
To put it in other way, the vertex with label $x=\l x_1, \ldots ,x_k\r\in \mathcal N$ is the $x_k$th children of the vertex with label $\l x_1, \ldots, x_{k-1}\r$, which is the $x_{k-1}$st children of its parent, and so on.  

Note that trees can be represented by the set of the labels of their vertices; the labels give all information about the edges. In the sequel we identify vertices with their labels, and trees with the set of labels. The set of finite rooted ordered trees is denoted by $\mathcal G$. We say that the vertex with label $x=\l x_1, \ldots, x_k\r\in \mathcal N$ belongs to the $k$th generation of $G \in\mathcal G$. The degree of a vertex $x$ in $G$ will be denoted by $\deg\l x, G\r$. 

\subsection{The random tree model with linear weight function}

We consider randomly growing trees. At each step we add a new vertex, which attaches to a randomly chosen, already existing vertex with one edge. The probability that a vertex of degree $d$ gets the new edge is proportional to a fixed linear function of $d$. 

To formulate this, let $\beta>-1$ be fixed, and let $G_n\in \mathcal G \ \l n\in \N\r$ be a sequence of random finite rooted ordered trees, such that $G_1=\lk \emptyset\rk$, and the following holds for all $n, k\in \N$, $x=\l x_1, \ldots, x_k\r\in G_n$ and $d=\deg\l x, G_n\r$.
\[\P\l \left.G_{n+1}=G_n\cup  \lk \l x_1, \ldots, x_k, d\r\rk\right\vert G_n\r= \frac{d+\beta}{S_n},\]
where $S_n=2n-2+n\beta=\sum_{v\in G_n} \l \deg\l v, G_n\r+\beta\r$. The condition $\beta>-1$ guarantees that the given probabilities are positive.

We say that the \textit{weight} of a vertex of degree $d$ is $d+\beta$. In other words, each vertex has weight $1+\beta$ when it is born, and its weight increases by 1 every time when it gives birth to a child.

\section{Main results}

Our goal is to describe the almost sure behavior and the limiting distribution of the degree of a fixed vertex $x$ in $G_n$ as $n\tv$. 

For the root, the moments of the limiting distribution are calculated in \cite{max}. It is proven that 
\begin{equation}\label{e1}\frac{\deg\l \emptyset, G_n\r}{n^{1/\l2+\beta\r}}\rightarrow \zeta_0\end{equation}
as $n\tv$, with probability 1, with a  positive random variable $\zeta_0$. Moreover, the following holds for every integer $k\geq 1$.
\begin{equation}\label{e10}\E \zeta_0^k=k!\cdot\frac{\Gamma\l1+ \frac{\beta}{2+\beta}\r}{\Gamma\l 1+\frac{k+\beta}{2+\beta}\r}\cdot\binom{k+\beta}{k}=
\frac{\Gamma\l1+ \frac{\beta}{2+\beta}\r}{\Gamma\l \beta+1\r}\cdot
\frac{\Gamma\l k+\beta+1\r}{\Gamma\l1+\frac{k+\beta}{2+\beta}\r}.\end{equation}

Our main result is the following. 

\begin{theorem}\label{t1}
 Let $k\in \Z_+$ and $x=\l x_1, \ldots, x_k\r\in \mathcal N$ be fixed. Then 
\[\frac{\deg\l x, G_n\r}{n^{1/\l2+\beta\r}}\rightarrow \zeta_x\]
with probability 1, with some positive random variable $\zeta_x$. The distribution of $\zeta_x$ is the same as the 
distribution of $\zeta_{0}\cdot \xi_1\cdot\ldots\cdot\xi_k$, where 
\begin{itemize}
 \item $\zeta_{0}, \xi_1, \ldots, \xi_k$ are independent random variables;
\item $\zeta_0$ is defined by equation $\l\ref{e1}\r$;
\item $\xi_1$ has distribution $\textit{Beta}\l 1+\beta, x_1-1\r$ if $x_1>1$; $\xi_1\equiv 1$ if $x_1=1$;
\item $\xi_s$ has distribution $\textit{Beta}\l 1+\beta, x_s\r$ for $2\leq s\leq k$.
\end{itemize}
\end{theorem}
\biz We prove the theorem by induction on $k$. 

For $k=1$, let $j\in \N$ be fixed. Vertex $x=j$ is the $j$th child of the root. 

For $j=1$ we have one edge, for symmetry reasons it is clear that $\zeta_1$ and $\zeta_0$ are identically distributed.

For $j>1$, assume that vertex $j$ appears in the $N$th step, that is, $j\in G_N\setminus G_{N-1}$. $N$ is a random positive integer. 
After the birth of vertex $j$, we divide the weight of the vertices into two parts, a "black" and  a "white" one. $W_n\l x\r$ and $B_n\l x\r$ denote the "white" and "black" weight of vertex $x$ in $G_n$, respectively, for $n\geq N$. The total weight of a vertex is equal to its actual degree plus $\beta$, thus we have
\begin{equation}\label{e11}\deg\l x, G_n\r+\beta=W_n\l x\r+B_n\l x\r\ \ \ \l n\geq N, x\in G_n\r.\end{equation}

We set 
\begin{equation*}
 \begin{split}
  W_N\l j\r&=1+\beta, \ \ \ B_N\l j\r=0;\\
W_N\l \emptyset\r&=1+\beta, \ \ \ B_N\l \emptyset\r=j-1;\\
W_N\l x\r&=0, \ \ \ \ \ \ \ \ B_N\l x\r=\deg\l x, G_N\r+\beta \ \ \ \l x \in G_N\setminus \lk \emptyset, j\rk\r. 
 \end{split}
\end{equation*}

This is possible, because vertex $j$ has weight $1+\beta$ when it is born, and the root has degree $j$ at the same time. All the other weights are colored black. 

Later on, when a new vertex appears, it gets black weight $1+\beta$, its total weight is black.  When an old vertex, $x$,  gets a new edge, its weight increases by 1. The color of this increment will be randomly chosen; the probability that the increment is white is the ratio of the white part to the total weight of $x$. We formulate this in the following way. For $n\geq N$, $k\in \N$, $x=\l x_1, \ldots, x_k\r\in G_n$ let $d=\deg\l x, G_n\r$ if $x\neq \emptyset$, and $d=1+\deg\l x, G_n\r$ if $x=\emptyset$. We set
\begin{equation}\label{e2}\begin{split}
   \P&\l \left.G_{n+1}=G_n\cup  \lk \l x_1, \ldots, x_k, d\r\rk, W_{n+1}\l x\r=W_n\l x\r+1\right\vert G_n\r\\&=\frac{d+\beta}{S_n}\cdot \frac{W_n\l x\r}{d+\beta}=\frac{W_n\l x\r}{S_n}, \\
\P&\l \left.G_{n+1}=G_n\cup  \lk \l x_1, \ldots, x_k, d\r\rk, B_{n+1}\l x\r=B_n\l x\r+1\right\vert G_n\r\\&=\frac{d+\beta}{S_n}\cdot \frac{B_n\l x\r}{d+\beta}=\frac{B_n\l x\r}{S_n}.
  \end{split}
\end{equation}

Of course, if $W_{n+1}\l x\r=W_n\l x\r+1$, then $B_{n+1}\l x\r=B_n\l x\r$, and if $B_{n+1}\l x\r=B_n\l x\r+1$, then $W_{n+1}\l x\r=W_n\l x\r$; the total weight is increased by 1 in both cases. 

Note that vertices, except the root and its $j$th child, never get white weight, which implies that
\[W_n\l x\r=0,\ \  B_n\l x\r=\deg\l x, G_n\r+\beta\ \ \ \l n\geq N, x\in G_n\setminus \lk \emptyset, j\rk\r. \]

On the other hand, the total weight of vertex $j$ is white when it is born, thus we have
\begin{equation}\label{e21}W_n\l j\r=\deg\l j, G_n\r+\beta, \ \ B_n\l j\r=0\ \ \ \l n\geq N\r.\end{equation}

We colored the weights in such a way that
\[W_N\l j\r=W_N\l \emptyset\r=1+\beta. \]

From formula $\l\ref{e2}\r$ it follows that the probability that the white weight of a vertex increases depends only on its actual white weight, independently  of the structure of $G_N$; and $S_n$ is deterministic. Thus, for symmetry reasons, we have  that if 
\[\frac{W_n\l\emptyset\r}{n^{1/\l 2+\beta\r}}\rightarrow\eta^0_j\]
almost surely for some positive random variable $\eta^0_j$, then
\[\frac{W_n\l j\r}{n^{1/\l 2+\beta\r}}\rightarrow\eta_j\]
almost surely as well, where $\eta^0_j$ and $\eta_j$ are identically distributed. Thus we will focus on the behavior of the root. 

Recall that $W_N\l \emptyset\r=1+\beta, B_N\l \emptyset\r=j-1$. From equation $\l \ref{e2}\r$ it follows that 
\begin{equation}\label{e3}\begin{split}
  \P&\l \left.W_{n+1}\l \emptyset\r=W_n\l \emptyset\r+1 \right\vert\deg\l\emptyset, G_{n+1}\r=\deg\l \emptyset, G_n\r+1\r\\&=\frac{W_n\l \emptyset\r}{W_n\l \emptyset\r+B_n\l\emptyset\r},\\
\P&\l \left.B_{n+1}\l \emptyset\r=B_n\l \emptyset\r+1 \right\vert\deg\l\emptyset, G_{n+1}\r=\deg\l \emptyset, G_n\r+1\r\\&=\frac{B_n\l \emptyset\r}{W_n\l \emptyset\r+B_n\l\emptyset\r}.
  \end{split}
\end{equation}

This corresponds to a P\'olya--Eggenbegger urn model. We have an urn with $a$ white and $b$ black balls. At each step, we draw a ball with uniform distribution, and put it back together with $c$ balls of the same color. It is well known (see e.g. \cite{Gouet}) that the proportion of the white balls converges almost surely to a $\textit{Beta}\l a/c, b/c\r$ distributed random variable. Moreover, saying that the number of the balls corresponds to "white" and "black" weights, this remains the same for positive, not necessarily integer weights. In this case, at each step, we choose a random color. The probability that white is chosen is equal to the actual proportion of the white weight in the urn. Then the weight of the selected color is increased by $c$. As equations $\l\ref{e3}\r$ show, this happens at the steps when the root gives birth to a child. Thus we can apply these results with $a=1+\beta$, $b=j-1$ and $c=1$, and we get that
\[\frac{W_n\l \emptyset\r}{W_n\l \emptyset\r+B_n\l \emptyset\r}\rightarrow \xi^0_j\] 
almost surely as $n\tv$, and $\xi^0_j$ has distribution $\textit{Beta}\l 1+ \beta, j-1\r$.

From this result, equations $\l \ref{e1}\r$, $\l \ref{e11}\r$ and the condition $\beta>-1$ it follows that 
\[\frac{W_n\l \emptyset\r}{n^{1/\l 2+\beta\r}}=\frac{W_n\l \emptyset\r}{W_n\l \emptyset\r+B_n\l \emptyset\r}\cdot\frac{\deg\l \emptyset, G_n\r+\beta}{n^{1/\l 2+\beta\r}}\rightarrow \xi^0_j\cdot\zeta_0=\eta^0_j\]
almost surely as $n\tv$. $\xi^0_j$ and $\zeta_0$ are independent, because the almost sure convergence of the proportion of the white weight does not depend on the behavior of the total weight of the root. As we have seen before, this implies that 
\[\frac{W_n\l j\r}{n^{1/\l 2+\beta\r}}\rightarrow \eta_j\]
almost surely as $n\tv$, where $\eta^0_j$ and $\eta_j$ are identically distributed. Thus, using equation $\l\ref{e21}\r$ we get that 
\[\frac{\deg\l j, G_n\r}{n^{1/\l 2+\beta\r}}=\frac{W_n\l j\r-\beta}{n^{1/\l 2+\beta\r}}\rightarrow\eta_j\]
almost surely as $n\tv$, and $\eta_j$ has the same distribution as $\zeta_0\cdot \xi_j$, where $\zeta_0$ is defined by equation $\l\ref{e1}\r$, $\xi_j$ has distribution $\textit{Beta}\l 1+\beta, j-1\r$, and finally, $\zeta_0$ and $\xi_j$ are independent. This completes the proof for the case $k=1$.

Assume that the statement is proven for some $k\geq 1$. We fix $x=\l x_1, \ldots, x_k\r\in \N$ and $x_{k+1}\in \mathbb Z_+$.  The induction step can be verified by a slight modification of the previous argument. The only difference is that $x$ has degree $x_{k+1}$ when its $x_{k+1}$st child is born, namely, one edge is attached to its parent, $\l x_1, \ldots, x_{k-1}\r$, and $x_{k+1}-1$ to its already existing children. This means that $a=1+\beta$ and $b=x_{k+1}$ in the urn model, and the last factor $\xi_{k+1}$ has distribution $\textit{Beta}\l 1+\beta, x_{k+1}\r$. \kesz

\section{A particular case and a variant}

\subsection{Albert--Barab\'asi tree}

For $\beta=0$ the probability that a given vertex of degree $d$ gets the new edge is proportional to $d$. 
This is the Albert--Barab\'asi random tree \cite{BA}. In this special case it is possible to determine the distribution of $\zeta_0$. Namely,  $2^{-1/2}\zeta_0$ is identically distributed with the absolute value of a standard normal random variable. To verify this, one can get the even moments of $2^{-1/2}\zeta_0$ from equation $\l\ref{e10}\r$, namely, 
\[\E\l \l 2^{-1/2}\zeta_0\r^{2k}\r=\frac{\Gamma\l 2k+1\r}{2^k\Gamma\l k+1\r}=\l 2k-1\r!!,\]
which are just the even moments of the standard normal distribution.

\subsection{Generalized PORT model}

Plane oriented random trees (PORT) are also rooted ordered trees. In this case the tree is embedded in the plane, and the left-to-right order of the children of the vertices is relevant (see e. g. \cite{Drmota}). The out-degree of vertex $v$ in tree $G$ is the number of its children, and it will be denoted by $\deg ^+\l v, G\r$. It is equal to the degree of the vertex minus one, except for the root, where it is the same as the degree. 

If a vertex $v$ has out-degree $d$, then there are $d+1$ possible ways to attach a new vertex to $v$. We get plane oriented random trees if all these possibilities are equally likely. Namely, the probability that a vertex of out-degree $d$ gets the new edge is proportional to $d+1$. In the generalized PORT model, this probability is proportional to $d+\beta$ for some $\beta>0$ \cite{Drmota1, Anna}.

We can get similar results to Theorem \ref{t1}.

\begin{theorem}\label{t3}
 Let $k\in \Z_+$ and $x=\l x_1, \ldots, x_k\r\in \mathcal N$ be fixed. Then 
\[\frac{\deg\l x, G_n\r}{n^{1/\l1+\beta\r}}\rightarrow \zeta_x\]
with probability 1, with some nondegerate random variable $\zeta_x$. The distribution of $\zeta_x$ is the same as the 
distribution of $\zeta_{0}\cdot \xi_1\cdot\ldots\cdot\xi_k$, where 
\begin{itemize}
 \item $\zeta_{0}, \xi_1, \ldots, \xi_k$ are independent random variables;
\item $\xi_s$ has distribution $\textit{Beta}\l \beta, x_s\r$ for $1\leq s\leq k$.
\end{itemize}
\end{theorem} 

After showing the existence of $\zeta_0$, the proof is a straightforward modification of the proof of Theorem $\ref{t1}$, therefore
we omit it.

To verify that $\zeta_0$ exists, we apply the results of Gouet about generalized P\'olya urn models \cite{Gouet93}. We color the root white and all 
the other vertices black. At the beginning, the root has weight $1+\beta$, while the only black vertex has weight $\beta$.
When the root gives birth to a child, its weight is increased by 1, and the black weight is increased by $\beta$. 
On the other hand, if the new vertex does not attach to the root, then the white weight does not change, while the 
black is increased by $1+\beta$. Thus we have the matrix 
\[R=\l 
\begin{array}{cc}
1 & \beta\\
0 & 1+\beta
       \end{array}
       \r.
\]
Recall that $\beta>0$. It is easy to check that all assumptions of Proposition 2.2 in \cite{Gouet93} hold, except that 
$\beta$ is not necessarily an integer. This is also easy to handle, and 
we obtain  for $W_n$, the white weight, after $n$ steps, that
\[\frac{W_n}{n^{1/\l 1+\beta\r}}\rightarrow Z\]
almost surely, where $Z$ is a nondegenerate random variable. This  guarantees the existence of $\zeta_0$. 

The moments of $\zeta_0$ can also be determined by the method applied in \cite{max}. Fix a positive integer $k$, $X_n$ denotes the degree of the root after $n$ steps, and let 
\[c_k=\prod_{i=1}^{n-1} \l 1+\frac{k}{i\l1+\beta\r-1}\r\asymp n^{k/\l 1+\beta\r}.\]
One can verify that 
\[Z_n=\frac{1}{c_n}\binom{X_n+k+\beta-1}{k} \ \ \ \l n=1, 2, \ldots\r\]
is a nonnegative martingale; it is convergent almost surely, and bounded in $L_p$ for every $p\geq 1$. Hence we obtain  by further calculations that 
\[\E\zeta_0^k=\frac{\Gamma\l\frac{\beta}{1+\beta}\r\Gamma\l k+\beta\r}{
\Gamma\l \frac{k+\beta}{1+\beta}\r\Gamma\l \beta\r}.\]

For $\beta=1$ we get the PORT model. In this special case $\zeta_0$ is identically distributed with $2\xi^{1/2}$, where $\xi$ has exponential distribution with expectation 1.

\end{document}